\documentclass[11pt]{article}
\usepackage{amssymb,latexsym,times}
\usepackage{mathptmx}
\usepackage{amscd}
\usepackage[all]{xy}

\newtheorem{Thm}{Theorem}[section]

\newtheorem{Prop}[Thm]{Proposition}
\newtheorem{Conj}[Thm]{Conjecture}
\newtheorem{example}[Thm]{Example}
\newtheorem{remark}[Thm]{Remark}
\newtheorem{Def}[Thm]{Definition}

\newcommand{\Z}{\mathbb{Z}}
\newcommand{\Q}{\mathbb{Q}}

\newcommand{\N}{\mathbb{N}}
\newcommand{\C}{\mathbb{C}}

\newcommand{\md}{\operatorname{mod}}

\usepackage{xcolor}
\usepackage{version}

\excludeversion{NB}

\newcommand{\bsm}{\begin{smallmatrix}}
\newcommand{\esm}{\end{smallmatrix}}

\newcommand{\bg}{\mathbf{g}}

\def\CC{{\mathcal C}}

\def\<{\langle\,}
\def\>{\,\rangle}

\def\eg{{\em e.g. }}

\def\g{\mathfrak g}

\def\<{\langle}
\def\>{\rangle}

\def\1{\mathbf 1}
\def\Gr{{\rm Gr}}

\def\a{\alpha}

\def\de{\delta}
\def\De{\Delta}
\def\Sl{\mathfrak{sl}}
\def\la{\lambda}
\def\AA{\mathcal{A}}

\def\1{{\mathbf 1}}
\def\ie{{\em i.e. }}
\def\hg{\widehat{\mathfrak{g}}}

\def\G{\Gamma}
\def\tG{\widetilde{\Gamma}}

\def\b{\beta}
\def\ga{\gamma}

\def\be{\mathbf e}
\def\bz{\mathbf z}

\def\hy{\widehat{y}}

\def\be{\mathbf{e}}
\def\by{\mathbf{y}}
\def\bv{\mathbf{v}}
\def\ba{\mathbf{a}}

\def\bdeux{2}
\def\bun{1}

\def\ds{\displaystyle}

\def\hatsl{\widehat{\mathfrak{sl}}}
\def\ev{\mathrm{ev}}
\def\hP{\widehat{P}}
\def\hM{\widehat{M}}
\def\hN{\widehat{N}}
\def\hmu{\widehat{\mu}}
\def\hla{\widehat{\la}}
\def\md{\mathrm{mod}}
\def\hby{\widehat{\mathbf{y}}}

\def\H{\mathcal{H}}
\def\LL{\mathcal{L}}

\begin{document}

\title{\bf Quantum affine algebras and cluster algebras}
\author{D. Hernandez\thanks{D. Hernandez is supported in part by the European Research Council under the European Union's 
Framework Programme H2020 with ERC Grant Agreement number 647353 QAffine.}, B. Leclerc}

\date{}

\maketitle

\begin{center}
{\it To Vyjayanthi Chari on her birthday}
\end{center}

\setcounter{tocdepth}{1}
{\small \tableofcontents}

\section{Introduction}

This article is an extended version of the minicourse given by the second author at the summer school of the conference 
\emph{Interactions of quantum affine algebras with cluster algebras, current algebras and categorification}, held in June 2018 in Washington.
The aim of the minicourse, consisting of three lectures, was to present a number of results and conjectures on certain monoidal categories of
finite-dimensional representations of quantum affine algebras, obtained by exploiting the fact that their Grothendieck rings have the natural
structure of a cluster algebra.

\section{A forerunner: Chari and Pressley's paper on $U_q(\hatsl_2)$}\label{SectSl2}

In \cite{CP1}, Chari and Pressley launched a systematic study of tensor categories of finite-dimen\-sio\-nal representations of quantum affine
algebras by investigating in detail the case of $U_q(\hatsl_2)$. They gave a classification of simple objects, as
well as a concrete description of them as tensor products of evaluation modules. They also gave a necessary and sufficient condition for 
such tensor products to be irreducible, and they described the composition factors of a reducible tensor product of two evaluation representations.

In retrospect, these results may be seen as providing a cluster algebra structure on the Grothen\-dieck ring of this category, 
predating by ten years the invention of cluster algebras by Fomin and Zelevinsky \cite{FZ1}.
We will therefore start our lectures by reviewing these results.  

\subsection{The Hopf algebra $U_q(\hatsl_2)$}

Throughout the paper, we fix $q\in \C^*$ not a root of unity.
The algebra $U_q(\hatsl_2)$ is generated over $\C$ by 
\[
 E_0,\ F_0,\ K_0,\ K_0^{-1},\ E_1,\ F_1,\ K_1,\ K_1^{-1},
\]
subject to the following relations:
\begin{eqnarray}
 &K_iK_i^{-1} = 1, \label{R1}\\
 &K_iK_j = K_jK_i, \\
 &K_iE_iK_i^{-1} = q^2E_i, \label{R2}\\
 &K_iE_jK_i^{-1} = q^{-2}E_j, \label{R3}\\
 &K_iF_iK_i^{-1} = q^{-2}F_i, \label{R4}\\
 &K_iF_jK_i^{-1} = q^2F_j, \label{R5}\\
 &E_iF_i-F_iE_i = \displaystyle\frac{K_i-K_i^{-1}}{q-q^{-1}},\label{R6}\\
 &E_iF_j-F_jE_i = 0,\\
 &E_i^3E_j-(q^2+1+q^{-2})E_i^2E_jE_i + (q^2+1+q^{-2})E_iE_jE_i^2 -E_jE_i^3 =0,\\
 &F_i^3F_j-(q^2+1+q^{-2})F_i^2F_jF_i + (q^2+1+q^{-2})F_iF_jF_i^2 -F_jF_i^3 =0,
\end{eqnarray}
where $i \not = j$ are indices in $\{0,1\}$. Moreover $U_q(\hatsl_2)$ is a Hopf algebra, with comultiplication $\De$
given by:
\begin{eqnarray*}
&&\De(E_i) = E_i \otimes K_i + 1 \otimes E_i,\\
&&\De(F_i) = F_i \otimes 1 + K_i^{-1} \otimes F_i,\\
&&\De(K_i) = K_i\otimes K_i.
\end{eqnarray*}
It follows that a tensor product of finite-dimensional $U_q(\hatsl_2)$-modules is again a $U_q(\hatsl_2)$-module. 

\subsection{Simple finite-dimensional $U_q(\hatsl_2)$-modules}

Let $E$, $F$, $K$, $K^{-1}$ denote the generators of $U_q(\Sl_2)$. (They are subject to the same relations as 
(\ref{R1})(\ref{R2})(\ref{R4})(\ref{R6}).) 

For every $a\in\C^*$, we have a surjective algebra homomorphism
$\mathrm{ev}_a : U_q(\hatsl_2) \to U_q(\Sl_2)$ such that:
\[
 \ev_a(E_1) = E,\quad \ev_a(F_1) = F,\quad \ev_a(E_0) = q^{-1}aF,\quad \ev_a(F_0) = qa^{-1}E.
\]
Hence, every simple finite-dimensional $U_q(\Sl_2)$-module $M$ becomes a finite-dimensional $U_q(\hatsl_2)$-module
$M(a)$ by pull-back through $\ev_a$.

It is well-known that the simple finite-dimensional $U_q(\Sl_2)$-modules\footnote{In these lectures, we will only consider type I representations
of quantum enveloping algebras. All representations can be obtained from the type I representations by twisting with some signs,
see \eg \cite[\S 10.1]{CP2}.} 
are classified by their dimension: for every $n\in \Z_{\ge 0}$ there is a unique (up to isomorphism) simple module $V_n$ with dimension $n+1$.
Therefore, pulling back by the evaluation morphisms $\ev_a$, we get for all $n\in \Z_{\ge 0}$ a one-parameter family of simple $U_q(\hatsl_2)$-modules 
$V_n(a)\ (a\in \C^*)$ with dimension $n+1$. The representations $V_0(a)$ are all equal to the trivial representation. Otherwise, for $n\ge 1$, 
the simple modules $V_n(a)$ and $V_n(b)$ are non-isomorphic if $a \not = b$.
The modules $V_n(a)$ are called \emph{evaluation modules}. 

\begin{Thm}{\rm\cite{CP1}}\label{Thm1} 
Every non-trivial simple finite-dimensional $U_q(\hatsl_2)$-module $M$ is isomorphic to a tensor product of evaluation modules, that is, 
\[
 M \simeq V_{n_1}(a_1)\otimes \cdots \otimes V_{n_k}(a_k)
\]
for some $k\in \Z_{> 0}$, $n_1,\ldots,n_k \in \Z_{>0}$, and $a_1,\ldots, a_k \in \C^*$.
\end{Thm}

Note that tensor products of evaluation modules are \emph{not} always irreducible. 
The next task is therefore to find some necessary and sufficient condition of irreducibility.
In order to formulate this condition we introduce the notion of a \emph{string}. 
This is a subset of $\C^*$ of the form:
\[
 \Sigma(n,a) := \{aq^{-n+1}, aq^{-n+3},\ldots, aq^{n-1}\},\qquad (n\in \Z_{\ge 0}, a\in \C^*).
\]
(In fact, $\Sigma(n,a)$ is nothing else than the set of roots of the Drinfeld polynomial of $V_n(a)$.)
We say that two strings $\Sigma_1$ and $\Sigma_2$ are \emph{in general position} if and only if
\begin{itemize}
 \item[(i)] $\Sigma_1 \cup \Sigma_2$ is not a string, or
 \item[(ii)] $\Sigma_1 \subseteq \Sigma_2$ or $\Sigma_2 \subseteq \Sigma_1$.
\end{itemize}

\begin{Thm}{\rm\cite{CP1}}\label{Thm2} 
The tensor product $V_{n_1}(a_1)\otimes \cdots \otimes V_{n_k}(a_k)$ is irreducible if and only if
for every $(i,j)\in \{1,\ldots ,k\}^2$ the strings $\Sigma(n_i,a_i)$ and $\Sigma(n_j,a_j)$ are in 
general position.
\end{Thm}
 
Two strings which are not in general position are called \emph{in special position}. What can we say
about the tensor product $V_{n_1}(a_1)\otimes V_{n_2}(a_2)$ when the strings $\Sigma_1:=\Sigma(n_1,a_1)$ and $\Sigma_2:=\Sigma(n_2,a_2)$ 
are in special position? It turns out that in this case the tensor product always has 
two non-isomorphic composition factors. These two irreducible modules are, by Theorem~\ref{Thm1} and Theorem~\ref{Thm2},
parametrized by two collections of strings in general position. Here is how to obtain them from $\Sigma_1$ and $\Sigma_2$.

Because of Theorem~\ref{Thm2}~(i), $\Sigma_3 := \Sigma_1\cup\Sigma_2$ is a string.
Clearly, $\Sigma_4 := \Sigma_1\cap\Sigma_2$ is also a string, contained in $\Sigma_3$.
Removing from $\Sigma_3$ the points of $\Sigma_4$ together with its two nearest neighbours, we are left with 
the union of two strings $\Sigma_5$ and $\Sigma_6$. 
It is easy to see that the two pairs of strings $(\Sigma_3,\Sigma_4)$ and $(\Sigma_5,\Sigma_6)$ are in general position.
For instance, if 
\[
\Sigma_1 = \{1,q^2,q^4,q^6,q^8\},\qquad \Sigma_2 = \{q^6,q^8,q^{10},q^{12},q^{14},q^{16}\} 
\]
then
\[
\Sigma_3 = \{1,q^2,q^4,q^6,q^8,q^{10},q^{12},q^{14},q^{16}\},\qquad \Sigma_4 = \{q^6,q^8\}  
\]
and 
\[
\Sigma_5 = \{1,q^2\}, \qquad \Sigma_6 = \{q^{12},q^{14},q^{16}\}. 
\]
We can then state:
\begin{Prop}{\rm\cite[Proposition 4.9]{CP1}}\label{Prop3}
Let $\Sigma_1$ and $\Sigma_2$ be two strings in special position.
With the above notation, in the Grothendieck ring the following relation holds:
\begin{equation}\label{Eq1}
  [V(\Sigma_1)\otimes V(\Sigma_2)] = [V(\Sigma_3)\otimes V(\Sigma_4)] + [V(\Sigma_5)\otimes V(\Sigma_6)].
\end{equation}
Here, $V(\Sigma_i)$ denotes the evaluation module whose associated string is $\Sigma_i$.
\end{Prop}

\subsection{Relation with cluster algebras}\label{ssecQACASl2}

A reader familiar with the definition of a cluster algebra will recognize in (\ref{Eq1}) 
an \emph{exchange relation}. Let us make this more precise. First note that if two strings are in special position, all their
points belong to the same class in $\C^*/q^{2\Z}$, that is, they all are of the form $aq^k$ for some fixed $a\in\C^*$ and
some $k\in 2\Z$. This motivates the following definition:
\begin{Def}{\rm\cite{HL1}}
Let $a\in \C^*$ and $\ell\in \Z_{>0}$. 
Let $\CC_{a,\ell}$ be the full subcategory of the category of finite-dimensional $U_q(\hatsl_2)$-modules whose objects $V$
satisfy:
\begin{quote}
 Every composition factor of $V$ is of the form $V_{n_1}(a_1)\otimes \cdots \otimes V_{n_k}(a_k)$ where all strings
 $\Sigma(n_i,a_i)$ are contained in $S:=\{a, aq^{-2}, \ldots, aq^{-2\ell}\}$.
\end{quote}
\end{Def}

The category $\CC_{a,\ell}$ depends only on $\ell$ up to isomorphism. We can therefore restrict ourselves to
the case $a=1$, and write $\CC_{1,\ell} = \CC_\ell$.
Then Theorem~\ref{Thm1}, Theorem~\ref{Thm2} and Proposition~\ref{Prop3} yield the following reformulation:
\begin{Thm}{\rm\cite{HL1}}
The category $\CC_\ell$ is a monoidal category, and its Grothendieck ring $K_0(\CC_\ell)$ has the structure of a cluster algebra
of finite type $A_\ell$ in the Fomin-Zelevinsky classification. More precisely, the cluster variables of $K_0(\CC_\ell)$
are the classes of the evaluation modules contained in $\CC_\ell$, the class $[V_{\ell + 1}(q^{-\ell})]$ being the only 
frozen variable. The cluster monomials are equal to the classes of the simple modules in $\CC_\ell$. Two cluster variables
are \emph{compatible} (\ie belong to the same cluster) if and only if the corresponding strings are in general position.
Otherwise they form an exchange pair with exchange relation given by (\ref{Eq1}).
\end{Thm}

\subsection{How can we generalize?}
In an attempt to extend these results from $U_q(\hatsl_2)$ to other quantum affine algebras, Chari and Pressley introduced
in \cite{CP4} the notion of a \emph{prime} module: this is a simple finite-dimensional module that cannot be factored as a
tensor product of modules of smaller dimension. It follows from Theorem~\ref{Thm1} that the prime $U_q(\hatsl_2)$-modules
are precisely the evaluation modules. This is no longer true for $U_q(\hatsl_3)$, and Chari and Pressley have constructed 
an infinite class of prime $U_q(\hatsl_3)$-modules which are not evaluation modules, see \cite{CP4}.  
In view of this, the following problems naturally arise. 

Let $\g$ be a simple Lie algebra over $\C$, and let $U_q(\hg)$ denote the corresponding untwisted quantum affine algebra.
\begin{itemize}
 \item[(P1)] What are the prime $U_q(\hg)$-modules?
 \item[(P2)] Which tensor products of prime $U_q(\hg)$-modules are simple?
\end{itemize}
It is known that Kirillov-Reshetikhin modules (see below Definition~\ref{DefKR}) are prime. 
The \emph{minimal affinization} modules introduced by Chari and Pressley \cite{CP3} as replacements 
for evaluation modules which do not exist outside type $A$, 
are also prime. But this is not a complete list as we shall see below.
Problem (P2) is also completely open.
In \cite{HL1} we proposed to use cluster algebras to shed new light on these questions.


\section{Reminder on finite-dimensional $U_q(\hg)$-modules}

\subsection{Cartan matrix}\label{ssecCartan}
Let $C=(c_{ij})_{i,j\in I}$ be the Cartan matrix of $\g$.
There is a diagonal matrix $D = \mbox{diag}(d_i\mid i\in I)$ with entries in $\Z_{>0}$ 
such that the product
\[
B=D\,C = (b_{ij})_{i,j\in I}
\]
is symmetric. We normalize $D$ so that $\min\{d_i \mid i\in I\} = 1$, and
we put $t:=\max\{d_i \mid i\in I\}$.
Thus 
\[
t =
\left\{
\begin{array}{cl}
1 & \mbox{if $C$ is of type $A_n$, $D_n$, $E_6$, $E_7$ or $E_8$}, \cr
         2 & \mbox{if $C$ is of type $B_n$, $C_n$ or $F_4$}, \cr
         3 & \mbox{if $C$ is of type $G_2$}.
\end{array}
\right.
\]
\begin{example}\label{example_Cartan}
{\rm
The Lie algebra $\g = \mathfrak{so}_7$, of type $B_3$ in the Cartan-Killing classification, has
Cartan matrix
\[
C = \pmatrix{
2 & -1 & 0 \cr
-1 & 2 & -1\cr
0 & -2 & 2}
\]
We have $D=\mbox{diag}(2,2,1)$ and the symmetric matrix $B$ is given by
\[
B = \pmatrix{
4 & -2 & 0 \cr
-2 & 4 & -2\cr
0 & -2 & 2}
\]
}
\end{example}

We denote by $\a_i\ (i\in I)$ the simple roots of $\g$, and by 
$\varpi_i\ (i\in I)$ the fundamental weights.
They are related by
\begin{equation}
\a_i = \sum_{j\in I} c_{ji} \varpi_j. 
\end{equation}

\subsection{Classification}
By Cartan-Killing theory, the simple finite-dimensional $\g$-modules are in one-to-one correspondence with their
highest weight, an element of the positive cone of integral dominant weights:
\[
 P_+ := \bigoplus_{i\in I} \N \varpi_i.
\]
We denote by $L(\lambda)$ the simple $\g$-module with highest weight $\la\in P_+$.

Chari and Pressley have obtained a similar classification of simple finite-dimensional $U_q(\hg)$-modules.
To formulate it, we introduce the cone 
\[
 \widehat{P}_+ := \bigoplus_{i\in I,\, a\in\C^*} \N (\varpi_i,a)
\]
of \emph{dominant loop-weights}.
\begin{Thm}{\rm\cite{CP2}}\label{Thm3} 
Up to isomorphism, the simple finite-dimensional $U_q(\hg)$-modules are in one-to-one correspondence with their
highest loop-weight, an element of $\widehat{P}_+$.
\end{Thm}
We denote by $L(\widehat{\lambda})$ the simple $U_q(\hg)$-module with highest loop-weight $\widehat{\la}\in \widehat{P}_+$.

\begin{example}\label{example_module}
{\rm
Let $\g = \mathfrak{so}_8$, of type $D_4$. Thus $I=\{1,2,3,4\}$, where we denote by 3 the trivalent node of the Dynkin diagram.
Then the $\g$-module $L(\varpi_3) = \bigwedge^2\C^8$ is of dimension 28. 

For $a \in \C^*$, the $U_q(\hg)$-module $L(\varpi_3,a)$ has dimension 29. This is a \emph{minimal affinization} of  
$L(\varpi_3)$ in the sense of \cite{CP3}.
}
\end{example}

\subsection{$q$-characters}
Finite-dimensional $\g$-modules $M$ are characterized by their character
\[
 \chi(M) := \sum_{\mu\in P} \dim(M_\mu) e^\mu,
\]
where $P := \oplus_{i\in I} \Z \varpi_i$ is the weight lattice, 
$M := \oplus_{\mu\in P} M_\mu$ is the weight space decomposition of $M$,
and $e^\mu$ is a formal exponential. So $\chi(M)$ is a Laurent polynomial
in the variables $y_i := e^{\varpi_i},\ (i\in I)$.

Similarly, finite-dimensional $U_q(\hg)$-modules $\widehat{M}$ have a loop-weight space
decomposition 
\[
 \widehat{M} := \oplus_{\hmu\in \hP} \widehat{M}_{\hmu}
\]
where $\hP:=\oplus_{i\in I,\, a\in\C^*} \Z (\varpi_i,a)$.
Frenkel-Reshetikhin introduced the \emph{$q$-character}
\[
 \chi_q(\hM) := \sum_{\hmu\in \hP} \dim(\hM_{\hmu}) e^{\hmu},
\]
a Laurent polynomial in the variables $Y_{i,a} := e^{(\varpi_i,a)},\ (i\in I,\ a\in\C^*)$.

\begin{Thm}{\rm\cite{FR}}\label{Thm4} 
Let $\hM$ and $\hN$ be two finite-dimensional $U_q(\hg)$-modules. 
The following are equivalent:
\begin{itemize}
 \item[(i)] $\chi_q(\hM) = \chi_q(\hN)$;
 \item[(ii)] $[\hM] = [\hN]$ in the Grothendieck ring $K_0(\md(U_q(\hg)))$;
 \item[(iii)] $\hM$ and $\hN$ have the same composition factors with the same multiplicities.
\end{itemize}
In particular, $q$-characters characterize \emph{simple} $U_q(\hg)$-modules up to isomorphism.
\end{Thm}

Note also that, because of Theorem~\ref{Thm4}, $\chi_q$ descends to an injective ring homomorphism from $K_0(\md(U_q(\hg)))$ to
the ring of Laurent polynomials $\Z[Y_{i,a}^{\pm1}\mid i\in I,\ a\in\C^*]$.

\begin{example}\label{example_qchar}
{\rm
Let $\g = \mathfrak{sl}_2$, of type $A_1$. Then, as is well known, we have
\begin{eqnarray*}
\chi(L(\varpi_1)) &=& \chi(V_1) \ =\  y_1 + y_1^{-1}\\
\chi(L(2\varpi_1)) &=& \chi(V_2) \ = \ y_1^2 + 1 + y_1^{-2}\\
\ldots && \ldots
\end{eqnarray*}
On the other hand, for any $a\in\C^*$,
\begin{eqnarray*}
\chi_q(L(\varpi_1,a)) &=& \chi_q(V_1(a)) \ =\  Y_{1,a} + Y_{1,aq^2}^{-1}\\
\chi_q(L(2(\varpi_1,a))) &=&  Y_{1,a}^2 + 2Y_{1,a}Y_{1,aq^2}^{-1} + Y_{1,aq^2}^{-2}\\
\chi_q(L((\varpi_1,a)+(\varpi_1,aq^2)) &=& \chi_q(V_2(aq)) \ =\  Y_{1,a}Y_{1,aq^2} + Y_{1,a}Y_{1,aq^4}^{-1} + Y_{1,aq^2}^{-1}Y_{1,aq^4}^{-1}
\end{eqnarray*}
This shows that $L((\varpi_1,a)+(\varpi_1,aq^2))$ is a minimal affinization of $L(2\varpi_i)$, but 
$L(2(\varpi_1,a))$ is not.
}
\end{example}

\begin{Def}\label{DefKR}
 For $i\in I$, $k\in\N$, $a\in\C^*$, set
 \[
  \hla_{k,a}^{(i)} := \sum_{j=0}^{k-1}(\varpi_i, aq^{2d_ij}) \in \hP_+.  
 \]
The simple $U_q(\hg)$-module $L\left(\hla_{k,a}^{(i)}\right)$ is called a \emph{Kirillov-Reshetikhin module}. 
We often write for short $W_{k,a}^{(i)} = L\left(\hla_{k,a}^{(i)}\right)$.
The modules $W_{1,a}^{(i)} = L((\varpi_i,a))$ are the \emph{fundamental} $U_q(\hg)$-modules.
\end{Def}

\subsection{$T$-systems}\label{ssect-T-syst} 

With the quantum affine algebra $U_q(\hg)$ is associated a system
of difference equations called a $T$-system \cite{KNS1}. 
Its unknowns are denoted by
\[
T^{(i)}_{k, r},\qquad (i \in I,\ k\in \N,\ r\in\Z). 
\]
We fix the initial boundary condition 
\begin{equation}
T^{(i)}_{0, r} = 1,\qquad (i\in I,\ r\in\Z).
\end{equation}
If $\g$ is of type $A_n, D_n, E_n$, the $T$-system
equations are
\begin{equation}
T^{(i)}_{k,r+1}T^{(i)}_{k,r-1} =  T^{(i)}_{k-1,r+1}T^{(i)}_{k+1,r-1} + 
\prod_{j :\ c_{ij}=-1} T^{(j)}_{k,r},\qquad (i \in I,\ k \ge 1,\ r\in\Z). 
\end{equation}
If $\g$ is not of simply laced type, the $T$-system equations are more complicated. 
They can be written in the form
\begin{equation}\label{Tsystem_general}
T^{(i)}_{k,r+d_i}T^{(i)}_{k,r-d_i} =  T^{(i)}_{k-1,r+d_i}T^{(i)}_{k+1,r-d_i} + 
S^{(i)}_{k,r},\qquad (i \in I,\ k \ge 1,\ r\in\Z), 
\end{equation}
where $S^{(i)}_{k,r}$ is defined as follows. 
If $d_i\ge 2$ then
\begin{equation}
 S^{(i)}_{k,r} = \prod_{j:\ c_{ji}=-1} T^{(j)}_{k,r} 
\prod_{j:\ c_{ji}\le -2} T^{(j)}_{d_ik,\ r-d_i+1}. 
\end{equation}
If $d_i=1$ and $t=2$, then
\begin{equation}
S^{(i)}_{k,r} = 
\left\{
\begin{array}{ll}
\ds\prod_{j:\ c_{ij}=-1} T^{(j)}_{k,r}
\prod_{j:\ c_{ij}= -2} T^{(j)}_{l,r} T^{(j)}_{l,r+2},
& \mbox{if $k=2l$,}\\[5mm]
\ds\prod_{j:\ c_{ij}=-1} T^{(j)}_{k,r}
\prod_{j:\ c_{ij}= -2} T^{(j)}_{l+1,r} T^{(j)}_{l,r+2}
& \mbox{if $k=2l+1$.}
\end{array}
\right.
\end{equation}
Finally, if $d_i=1$ and $t=3$, that is, if $\g$ is of type $G_2$,
denoting by $j$ the other vertex we have $d_j=3$ and
\begin{equation}
S^{(i)}_{k,r} = 
\left\{
\begin{array}{ll}
T^{(j)}_{l,r} T^{(j)}_{l,r+2} T^{(j)}_{l,r+4} & \mbox{if $k=3l$,}\\[2mm]
T^{(j)}_{l+1,r} T^{(j)}_{l,r+2} T^{(j)}_{l,r+4} & \mbox{if $k=3l+1$,}\\[2mm]
T^{(j)}_{l+1,r} T^{(j)}_{l+1,r+2} T^{(j)}_{l,r+4} & \mbox{if $k=3l+2$.}
\end{array}
\right.
\end{equation}

\begin{example}\label{example_Tsystem}
{\rm
Let $\g$ be of type $B_2$. 
The Cartan matrix is 
\[
 C =
\pmatrix{2 & -1\cr
-2 & 2}
\]
and we have $d_1=2$ and $d_2=1$.
The $T$-system reads:
\[
\begin{array}{lcll}
T^{(\bun)}_{k,r+2}T^{(\bun)}_{k,r-2} &=&  T^{(\bun)}_{k-1,r+2}T^{(\bun)}_{k+1,r-2} + 
T^{(\bdeux)}_{2k,r-1},& (k \ge 1,\ r\in\Z),\\[2mm] 
T^{(\bdeux)}_{2l,r+1}T^{(\bdeux)}_{2l,r-1} &=&  T^{(\bdeux)}_{2l-1,r+1}T^{(\bdeux)}_{2l+1,r-1} + 
T^{(\bun)}_{l,r} T^{(\bun)}_{l,r+2}, &(l \ge 1,\ r\in\Z), \\[2mm]
T^{(\bdeux)}_{2l+1,r+1}T^{(\bdeux)}_{2l+1,r-1} &=&  T^{(\bdeux)}_{2l,r+1}T^{(\bdeux)}_{2l+2,r-1} + 
T^{(\bun)}_{l+1,r} T^{(\bun)}_{l,r+2}, &(l \ge 0,\ r\in\Z). 
\end{array}
\]

} 
\end{example}

It was conjectured in \cite{KNS1}, and proved in \cite{N2} (for $\g$ of type $A,D,E$) and \cite{H} (general case),
that the $q$-characters of the Kirillov-Reshetikhin modules of $U_q(\hg)$ satisfy the 
corresponding $T$-system. More precisely, we have
\begin{Thm}[\cite{N1}\cite{H}]\label{thm_Tsystem}
For $i \in I,\ k\in \N,\ r\in\Z$,
\[
T^{(i)}_{k,r} = \chi_q\left(W^{(i)}_{k,q^r}\right),
\]
is a solution of the $T$-system in the ring $\Z\left[Y_{i,q^r}^{\pm1}\mid (i,r) \in I\times\Z\right]$.
Equivalently, by Theorem~\ref{Thm4}, 
\[
T^{(i)}_{k,r} = \left[W^{(i)}_{k,q^r}\right], 
\]
is a solution of the $T$-system in the Grothendieck ring $K_0(\md(U_q(\hg)))$.
\end{Thm}

\begin{remark} 
{\rm
(i)\ For $\g = \Sl_2$, Theorem~\ref{thm_Tsystem} is a particular case of Proposition~\ref{Prop3}.

\medskip\noindent
(ii)\ Theorem~\ref{thm_Tsystem} allows to calculate inductively $q$-characters of Kirillov-Reshetikhin modules
in terms of $q$-characters of fundamental modules.
Note, however, that it is not straightforward to compute the $q$-characters of the fundamental modules in type $E_8$ or $F_4$, say.
An algorithm has been obtained by Frenkel and Mukhin \cite{FM}. Another one, based on cluster mutation, is described in \cite{HL2}.

\medskip\noindent
(iii)\ The $T$-system in the Grothendieck ring comes from a non-split short exact sequence
\[
0\rightarrow \mathcal{S}^{(i)}_{k,r} \rightarrow W^{(i)}_{k,q^{r - d_i}} \otimes W^{(i)}_{k,q^{r + d_i}}\rightarrow 
W^{(i)}_{k-1,q^{r + d_i}}\otimes W^{(i)}_{k+1,q^{r + d_i}} \rightarrow 0,
\]
where the module $\mathcal{S}^{(i)}_{k,r}$ is defined as the tensor product of Kirillov-Reshetikhin modules 
associated with the factors of $S^{(i)}_{k,r}$ above (it does not depend on the order of the tensor product up to isomorphism). 
The representation $W^{(i)}_{k,q^{r - d_i}} \otimes W^{(i)}_{k,q^{r + d_i}}$ is of length $2$. 
By a general result of Chari on tensor products of Kirillov-Reshetikhin modules \cite{C}, it is cyclic generated by the tensor product 
of the highest weight vectors and so it is indecomposable.}\end{remark}


\section{Quivers, subcategories, and cluster algebras}

Following \cite{HL2}, we attach an infinite quiver to $U_q(\hg)$, and we define some subcategories of the category 
of finite-dimensional $U_q(\hg)$-modules. We then introduce cluster algebras corresponding to finite segments of 
this infinite quiver.

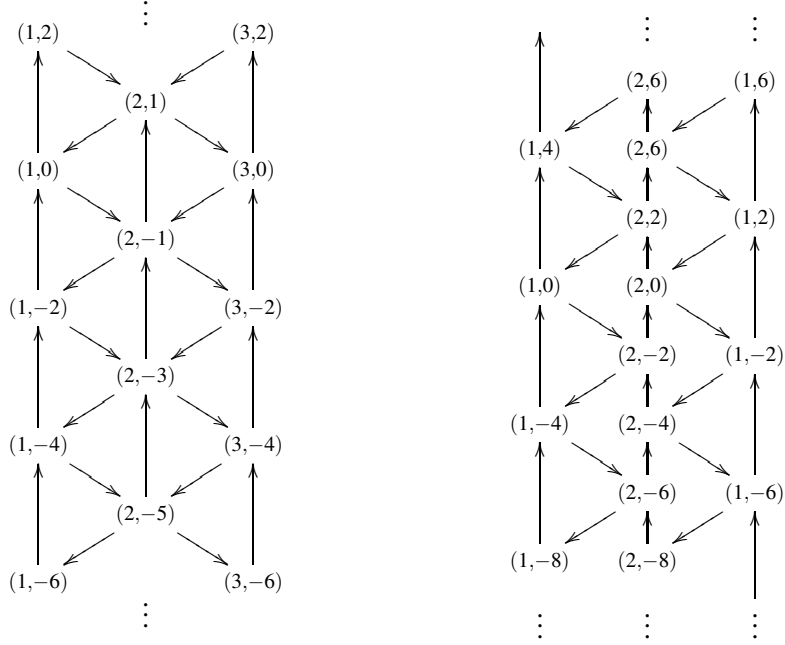
\begin{figure}[t]
\[
\def\objectstyle{\scriptstyle}
\def\lablestyle{\scriptstyle}
\xymatrix@-1.0pc{
&&&&\\
&{(1,2)}\ar[rd]&{}\save[]+<0cm,2ex>*{\vdots}\restore
&\ar[ld] (3,2) 
\\
&&\ar[ld] (2,1) \ar[rd]&&
\\
&\ar[uu]{(1,0)}\ar[rd]&
&\ar[ld] (3,0) \ar[uu]
\\
&&\ar[uu]\ar[ld] (2,-1) \ar[rd]&&
\\
&\ar[uu](1,-2) \ar[rd] &&\ar[ld] (3,-2)\ar[uu]
\\
&&\ar[ld] \ar[uu](2,-3) \ar[rd]&&
\\
&\ar[uu](1,-4) \ar[rd] &&\ar[ld] (3,-4)\ar[uu]
\\
&&\ar[ld] \ar[uu](2,-5) \ar[rd]&&
\\
&\ar[uu](1,-6) &{}\save[]+<0cm,-2ex>*{\vdots}\restore&\ar[uu] (3,-6) 
\\
}
\qquad\qquad
\xymatrix@-1.0pc{
&&&&\\
&&{}\save[]+<0cm,0ex>*{\vdots}\restore
& {}\save[]+<0cm,0ex>*{\vdots}\restore
\\
&&\ar[ld] (\bdeux,6)& \ar[ld](\bun,6)&
\\
&{(\bun,4)}\ar[rd]\ar[uu] &\ar[u] (\bdeux,6) \ar[rd]&&
\\
&&\ar[u]\ar[ld] (\bdeux,2) &\ar[ld] (\bun,2) \ar[uu]  &
\\
&\ar[uu](\bun,0)\ar[rd]   & \ar[u](\bdeux,0) \ar[rd]&&
\\
&&\ar[u]\ar[ld] (\bdeux,-2) &\ar[ld] (\bun,-2) \ar[uu] &
\\
&\ar[uu]\ar[rd](\bun,-4) & \ar[u]\ar[rd](\bdeux,-4)&&
\\
&&\ar[u]\ar[ld] (\bdeux,-6) &\ar[ld] (\bun,-6) \ar[uu]  &
\\
&\ar[uu](\bun,-8)   & \ar[u](\bdeux,-8) &&
\\
&{}\save[]+<0cm,0ex>*{\vdots}\restore  &{}\save[]+<0cm,0ex>*{\vdots}\restore&{}\save[]+<0cm,0ex>*{\vdots}\ar[uu]\restore 
\\
}
\]
\caption{\label{Fig0} {\it The quivers $\G$ in type $A_3$ and $B_2$.}}
\end{figure}

\subsection{Quivers}\label{ssectQuivers}
Put $\widetilde{V} = I \times \Z$.
We introduce a quiver $\tG$ with vertex set $\widetilde{V}$.
Recall the symmetric matrix $B= (b_{ij})_{i,j\in I}$ of \S\ref{ssecCartan}. The arrows of $\tG$ are given by
\[
((i,r) \to (j,s)) 
\quad \Longleftrightarrow \quad
(b_{ij}\not = 0 
\quad \mbox{and} \quad
s=r+b_{ij}).
\]
It is easy to check that the oriented graph $\tG$ has two isomorphic connected components. 
We pick one of them and call it $\Gamma$. The vertex set of $\Gamma$ is denoted by $V$.
Examples in type $A_3$ and $B_2$ are shown in Figure~\ref{Fig0}.

\subsection{Subcategories}

First, using the vertex set $V$, we introduce
\[
\hP_{+,\Z} := \bigoplus_{(i,r)\in V} \N (\varpi_i,q^{r+d_i}),
\]
a discrete subset of the positive cone $\hP_+$ of loop-weights.
\begin{Def}\label{defCZ}
Let $\CC_\Z$ be the full subcategory of the category of finite-dimensional $U_q(\hg)$-modules whose objects $M$ satisfy:
\begin{quote}
Every composition factor of $M$ is of the form $L(\hla)$ with $\hla\in\hP_{+,\Z}$.
\end{quote} 
\end{Def}
By \cite{HL2}, $\CC_\Z$ is a monoidal subcategory, \ie it is stable under tensor products.
Moreover, it is known that every simple finite-dimensional $U_q(\hg)$-module can be written
as a tensor product of simple objects of $\CC_\Z$ with spectral shifts. Therefore it is enough
to study the simple objects of $\CC_\Z$.

In order to relate $\CC_\Z$ with cluster algebras of finite rank, we need to ``truncate'' it, as we did in 
\S\ref{ssecQACASl2} for $\g = \Sl_2$.
Fix $\ell \in \N$, and put 
\[
\hP_{+,\ell} := \bigoplus_{(i,r)\in V,\, -2\ell-1\le r+d_i\le 0} \N (\varpi_i,q^{r+d_i}).
\]
\begin{Def}
Let $\CC_\ell$ be the full subcategory of the category of finite-dimensional $U_q(\hg)$-modules whose objects $M$ satisfy:
\begin{quote}
Every composition factor of $M$ is of the form $L(\hla)$ with $\hla\in\hP_{+,\ell}$.
\end{quote} 
\end{Def}
Again by \cite{HL2}, $\CC_\ell$ is a monoidal category and its Grothendieck ring is a polynomial ring in finitely many variables,
namely, the classes of the fundamental modules contained in $\CC_\ell$:
\[
 K_0(\CC_\ell) = \Z\left[[L\left((\varpi_i,q^{r+d_i})\right)]\mid (i,r)\in V,\ -2\ell-1\le r+d_i\le 0\right].
\]

\subsection{Cluster algebras}

We refer the reader to \cite{FZsurvey} and \cite{GSV} for an introduction to cluster algebras, and for any un\-defined terminology.

Let $\Gamma_\ell$ denote the full subquiver of $\Gamma$ with vertex set 
\[
 V_\ell := \{(i,r)\in V \mid -2\ell-1\le r+d_i\le 0\}. 
\]
Let $\bz_\ell := \{z_{(i,r)}\mid (i,r) \in V_\ell\}$ be a set of commuting indeterminates indexed by $V_\ell$. The pair
$(\bz_\ell, \Gamma_\ell)$ can be regarded as a \emph{seed}, in the sense of \cite{FZsurvey}.
\begin{Def}
Let $\AA_\ell\subset \Q(\bz_\ell)$ be the cluster algebra with initial seed $(\bz_\ell, \Gamma_\ell)$, where we consider
the variables $z_{i,r}$ with $r-d_i < -2\ell - 1$ as \emph{frozen variables}.
\end{Def}


\section{Main conjecture}

\subsection{Statements and examples}
The category $\CC_\ell$ and the cluster algebra $\AA_\ell$ are related as follows:

\begin{Thm}{\rm\cite{HL2}}\label{Thm5}
For $(i,r)\in V_\ell$ put $m_{i,r} := \max\{k\mid r+(2k+1)d_i\le 0\} + 1$.
The assignment $z_{i,r} \mapsto \left[W_{m_{i,r}, q^{r+d_i}}^{(i)}\right]$ extends to a
ring isomorphism $\iota_\ell : \AA_\ell \to K_0(\CC_\ell)$. 
\end{Thm}

We can now formulate our main conjecture from \cite{HL1, HL2}.

\begin{Conj}\label{mainConj}
\begin{itemize}
 \item[(i)] The isomorphism $\iota_\ell$ maps the subset of cluster monomials of $\AA_\ell$ into
the subset of classes of simple objects of $\CC_\ell$.
 \item[(ii)] The isomorphism $\iota_\ell$ maps the subset of cluster variables of $\AA_\ell$ into
the subset of classes of \emph{prime} simple objects of $\CC_\ell$.
\end{itemize}
\end{Conj}

In the situation of Theorem~\ref{Thm5} and Conjecture~\ref{mainConj}, we say that $\CC_\ell$ is 
a \emph{monoidal categorification} of the cluster algebra $\AA_\ell$ \cite{HL1}.
We illustrate the conjecture with simple examples.

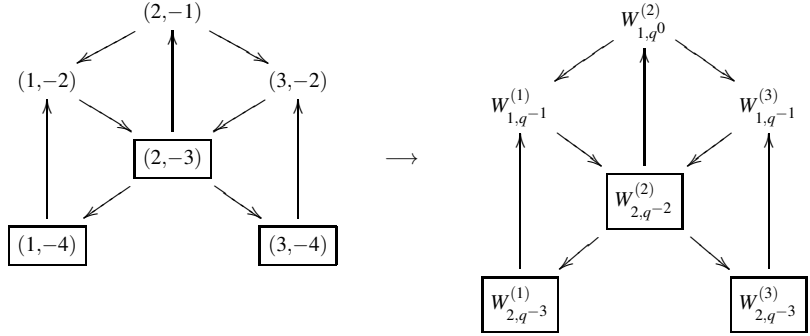
\begin{figure}[t]
\[
\def\objectstyle{\scriptstyle}
\def\lablestyle{\scriptstyle}
\xymatrix@-1.0pc{
&&&&&\\
&&\ar[ld] (2,-1) \ar[rd]&&
\\
&(1,-2) \ar[rd] &&\ar[ld] (3,-2)
\\
&&\ar[ld] \ar[uu]\fbox{$\scriptstyle(2,-3)$} \ar[rd]&&\longrightarrow
\\
&\ar[uu]\fbox{$\scriptstyle(1,-4)$} && \fbox{$\scriptstyle(3,-4)$}\ar[uu]
}  
\xymatrix@-1.0pc{
&&&\\
&\ar[ld] W_{1,q^0}^{(2)} \ar[rd]&&
\\
W_{1,q^{-1}}^{(1)} \ar[rd] &&\ar[ld] W_{1,q^{-1}}^{(3)}
\\
&\ar[ld] \ar[uu]\fbox{$\scriptstyle W_{2,q^{-2}}^{(2)}$} \ar[rd]&&
\\
\ar[uu]\fbox{$\scriptstyle W_{2,q^{-3}}^{(1)}$} && \fbox{$\scriptstyle W_{2,q^{-3}}^{(3)}$}\ar[uu]
}
\]
\caption{\label{Fig1} {\it The isomorphism $\iota_\ell$ in type $A_3$ for $\ell = 1$.}}
\end{figure}

\begin{example}\label{example_conjecture_A3_l1}
{\rm
We take $\g = \Sl_4$, of type $A_3$, and we choose $\ell = 1$. The isomorphism $\iota_1$ is 
displayed in Figure~\ref{Fig1}, which shows the image of the initial seed of $\AA_1$. 
For instance, $\iota_1(z_{2,-1}) = \left[W_{1,q^0}^{(2)}\right]$.
The 3 frozen variables are marked with a box.

The cluster variable $z_{2,-1}^*$ obtained by mutating $z_{2,-1}$ in this initial seed is
given by the exchange relation:
\[
z_{2,-1}z_{2,-1}^*  = z_{2,-3} + z_{1,-2}z_{3,-2}, 
\]
which translates under $\iota_1$ into the $T$-system equation:
\[
\left[W_{1,q^0}^{(2)}\right] \left[W_{1,q^{-2}}^{(2)}\right]
= \left[W_{2,q^{-2}}^{(2)}\right] + \left[W_{1,q^{-1}}^{(1)} \right] \left[W_{1,q^{-1}}^{(3)} \right].
\]
Thus, $\iota_1(z_{2,-1}^*) = \left[W_{1,q^{-2}}^{(2)}\right]$.

In this case, $\AA_1$ is a cluster algebra of finite type $A_3$ in the Fomin-Zelevinsky classification,
and Conjecture~\ref{mainConj} is proved \cite{HL1}. Moreover, since $\AA_1$ has finite type, the inclusions
of Conjecture~\ref{mainConj} are in fact bijections.
The algebra $\AA_1$ has 9 cluster variables plus 3 frozen variables. 
The prime simple modules of $\CC_1$ corresponding to cluster variables are
\[
\begin{array}{lll}
W_{1,q^{-1}}^{(1)} = L((\varpi_1,q^{-1})),&
W_{1,q^{-2}}^{(2)} = L((\varpi_2,q^{-2})),&
W_{1,q^{-1}}^{(3)} = L((\varpi_3,q^{-1})),\\[3mm]
W_{1,q^{-3}}^{(1)} = L((\varpi_1,q^{-3})),&
W_{1,q^{0}}^{(2)} = L((\varpi_2,q^{0})),&
W_{1,q^{-3}}^{(3)} = L((\varpi_3,q^{-3})),\\[3mm]
L((\varpi_1,q^{-3})+(\varpi_2,q^{0})), &
L((\varpi_1,q^{-3})+(\varpi_2,q^{0})+(\varpi_3,q^{-3})),&
L((\varpi_2,q^{0})+ (\varpi_3,q^{-3})).
\end{array}
\]
There are 6 fundamental modules, and 2 minimal affinizations, but
the 70-dimensional module $L((\varpi_1,q^{-3})+(\varpi_2,q^{0})+(\varpi_3,q^{-3}))$, which restricts to
$L(\varpi_1+\varpi_2+\varpi_3) \oplus L(\varpi_2)$ as a $U_q(\Sl_4)$-module, is \emph{not}
a minimal affinization.

\begin{figure}[t]
\begin{center}
\setlength{\unitlength}{1.70pt} 
\begin{picture}(90,110)(0,0) 
\thicklines

\put(0,0){\line(1,0){60}} 
\put(0,0){\line(0,1){40}} 
\put(60,0){\line(0,1){20}} 
\put(60,0){\line(1,1){30}} 
\put(0,40){\line(1,0){40}} 
\put(0,40){\line(1,3){20}} 
\put(60,20){\line(-1,1){20}} 
\put(60,20){\line(1,3){10}} 
\put(90,30){\line(0,1){40}} 
\put(40,40){\line(1,3){10}} 
\put(70,50){\line(1,1){20}} 
\put(70,50){\line(-1,1){20}} 
\put(50,70){\line(-1,3){10}} 
\put(90,70){\line(-1,1){40}} 
\put(20,100){\line(1,0){20}} 
\put(20,100){\line(1,1){10}} 
\put(30,110){\line(1,0){20}} 
\put(40,100){\line(1,1){10}} 
 
\thinlines 
\multiput(0,0)(1.5,1.5){20}{\circle*{0.5}} 
\multiput(30,30)(2,0){30}{\circle*{0.5}} 
\multiput(30,30)(0,2){40}{\circle*{0.5}} 
 
\put(35,105){\makebox(0,0){$\alpha_2$}} 
\put(28,70){\makebox(0,0){$\alpha_1+\alpha_2$}} 
\put(63,80){\makebox(0,0){$\alpha_2+\alpha_3$}} 
\put(55,45){\makebox(0,0){$\scriptstyle \alpha_1+\alpha_2+\alpha_3$}} 
\put(32,18){\makebox(0,0){$\alpha_1$}} 
\put(77,37){\makebox(0,0){$\alpha_3$}}

\put(0,0){\circle*{2}} 
\put(60,0){\circle*{2}} 
\put(60,20){\circle*{2}} 
\put(30,30){\circle*{2}} 
\put(90,30){\circle*{2}} 
\put(0,40){\circle*{2}} 
\put(40,40){\circle*{2}} 
\put(70,50){\circle*{2}} 
\put(50,70){\circle*{2}} 
\put(90,70){\circle*{2}} 
\put(20,100){\circle*{2}} 
\put(40,100){\circle*{2}} 
\put(30,110){\circle*{2}} 
\put(50,110){\circle*{2}}

\end{picture}  
\end{center}
\caption{\label{Fig2} {\it The associahedron in type $A_3$.}}
\end{figure}
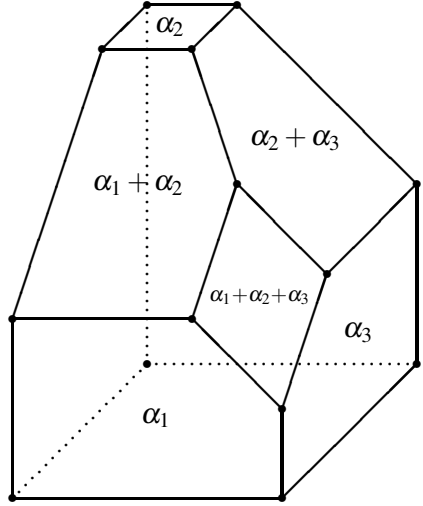

By \cite{FZ2}, there is a bijection between the set of cluster variables of $\AA_1$ and the set of almost positive roots of a
root system of type $A_3$. A natural bijection was given in \cite{HL1}:
\[
\begin{array}{lll}
L((\varpi_1,q^{-1})) \leftrightarrow -\alpha_1&
L((\varpi_2,q^{-2})) \leftrightarrow -\alpha_2&
L((\varpi_3,q^{-1}))\leftrightarrow -\alpha_3\\[3mm]
L((\varpi_1,q^{-3})) \leftrightarrow \alpha_1&
L((\varpi_2,q^{0})) \leftrightarrow \alpha_2&
L((\varpi_3,q^{-3})) \leftrightarrow \alpha_3\\[3mm]
L((\varpi_1,q^{-3})+(\varpi_2,q^{0})) \leftrightarrow \alpha_1 + \alpha_2 &
L((\varpi_2,q^{0})+ (\varpi_3,q^{-3})) \leftrightarrow \alpha_2+\alpha_3&\\[3mm]
L((\varpi_1,q^{-3})+(\varpi_2,q^{0})+(\varpi_3,q^{-3}))& \leftrightarrow \alpha_1+\alpha_2 + \alpha_3&
\end{array} 
\]
Following \cite{FZ2} and using this bijection one can read the 14 clusters of $\AA_1\cong K_0(\CC_1)$ from the associahedron of Figure~\ref{Fig2}.
Every cluster variable corresponds to a face, indicated by the attached almost positive root (the negative simple root $-\alpha_i$ labels
the unique rear face parallel to the face labelled by $\alpha_i$). Every cluster corresponds to a vertex and consists of the 3
faces adjacent to it. For instance, there is a cluster 
\[
 \{\alpha_1,\alpha_3,-\alpha_2\} \equiv \{ L((\varpi_1,q^{-3})),\, L((\varpi_3,q^{-3})),\, L((\varpi_2,q^{-2}))\}.
\]
The neat final result is that every simple module of $\CC_1$ is a tensor product of prime simple modules belonging to a single cluster,
and of frozen simple modules (corresponding to the frozen variables of $\AA_1$).
}
\end{example}

\begin{figure}[t]
\[
\def\objectstyle{\scriptstyle}
\def\lablestyle{\scriptstyle}
\xymatrix@-1.0pc{
&&&&\\
&&\ar[ld] (\bdeux,-2) &\ar[ld] (\bun,-2)  &
\\
&\ar[rd]\fbox{$\scriptstyle(\bun,-4)$} & \ar[u]\ar[rd](\bdeux,-4)&&\longrightarrow
\\
&&\ar[u] \fbox{$\scriptstyle(\bdeux,-6)$} & \fbox{$\scriptstyle(\bun,-6)$} \ar[uu]  &
}
\xymatrix@-1.0pc{
&&&&\\
&&\ar[ld] W_{1,q^{-1}}^{(2)} &\ar[ld] W_{1,q^{0}}^{(1)} &
\\
&\ar[rd]\fbox{$\scriptstyle W_{1,q^{-2}}^{(1)}$} & \ar[u]\ar[rd]W_{2,q^{-3}}^{(2)}&&
\\
&&\ar[u] \fbox{$\scriptstyle W_{3,q^{-5}}^{(2)}$} & \fbox{$\scriptstyle W_{2,q^{-4}}^{(1)}$} \ar[uu]  &
}
\]
\caption{\label{Fig3} {\it The isomorphism $\iota_\ell$ in type $B_2$ for $\ell = 2$.}}
\end{figure}
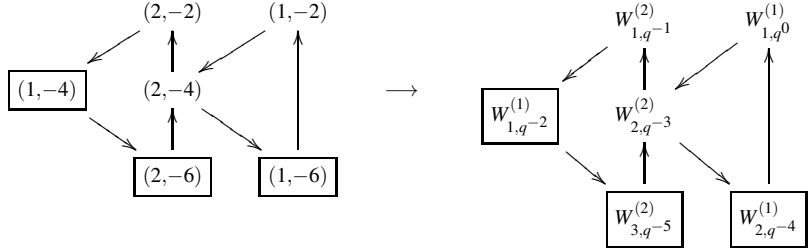

\begin{example}\label{example_conjecture_B2_l1}
{\rm
We take $\g = \mathfrak{so}_5$, of type $B_2$, and we choose $\ell = 2$.
The isomorphism $\iota_2$ is 
displayed in Figure~\ref{Fig3}, which shows the image of the initial seed of $\AA_2$. 
For instance, $\iota_2(z_{2,-2}) = \left[W_{1,q^{-1}}^{(2)}\right]$.
The 3 frozen variables are marked with a box.

Again, $\AA_2$ has finite cluster type $A_3$. So $\CC_2$ has 12 prime objects, namely, the 6 Kirillov-Reshetikhin
modules of the initial seed, together with
\[
 W_{1,q^{-4}}^{(1)},\
 W_{1,q^{-3}}^{(2)},\
 W_{1,q^{-5}}^{(2)},\
 W_{2,q^{-5}}^{(2)},\
 L((\varpi_1,q^0)+(\varpi_2,q^{-5})),\
 L((\varpi_1,q^0)+(\varpi_2,q^{-3})+(\varpi_2,q^{-5})).
\]
A bijection between the highest loop-weights of the unfrozen primes and the almost positive roots of type $A_3$, 
allowing to determine the clusters using the associahedron
as in Example~\ref{example_conjecture_A3_l1}, is for instance:
\[
\begin{array}{lll}
(\varpi_1,q^{-4}) \leftrightarrow -\alpha_1&
(\varpi_2,q^{-3})+(\varpi_2,q^{-1}) \leftrightarrow -\alpha_2&
(\varpi_2,q^{-1})\leftrightarrow -\alpha_3\\[3mm]
(\varpi_1,q^{0}) \leftrightarrow \alpha_1&
(\varpi_2,q^{-5}) \leftrightarrow \alpha_2&
(\varpi_2,q^{-3}) \leftrightarrow \alpha_3\\[3mm]
(\varpi_1,q^{0})+(\varpi_2,q^{-5}) \leftrightarrow \alpha_1 + \alpha_2 &
(\varpi_2,q^{-3})+ (\varpi_2,q^{-5}) \leftrightarrow \alpha_2+\alpha_3&\\[3mm]
(\varpi_1,q^{-0})+(\varpi_2,q^{-3})+(\varpi_2,q^{-5})& \leftrightarrow \alpha_1+\alpha_2 + \alpha_3&
\end{array}
\]

}
\end{example}


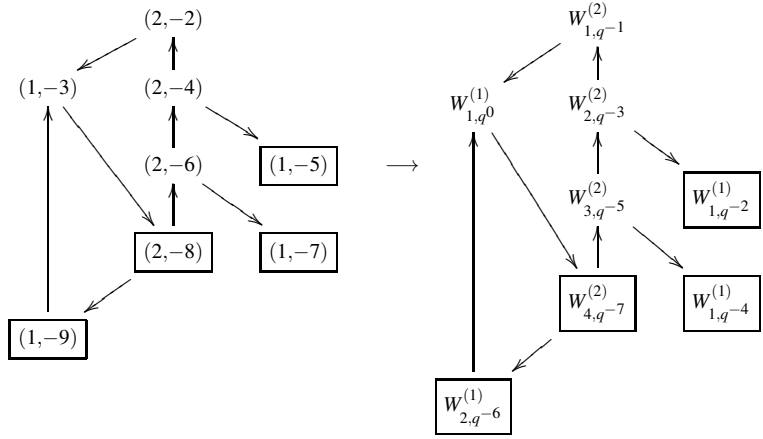
\begin{figure}[t]
\[
\def\objectstyle{\scriptstyle}
\def\lablestyle{\scriptstyle}
\xymatrix@-1.0pc{
&\ar[ld] (\bdeux,-2)& &
\\
{(\bun,-3)}\ar[ddr] &\ar[u] (\bdeux,-4) \ar[rd]&&
\\
&\ar[u] \ar[rd](\bdeux,-6) & \fbox{$\scriptstyle{(\bun,-5)}$} &\longrightarrow
\\
  &\fbox{$\scriptstyle{(\bdeux,-8)}$} \ar[u]\ar[dl]& \fbox{$\scriptstyle{(\bun,-7)}$} &
\\
\fbox{$\scriptstyle{(\bun,-9)}$} \ar[uuu]&&&
}
\xymatrix@-1.0pc{
&\ar[ld] W_{1,q^{-1}}^{(2)}& 
\\
W_{1,q^0}^{(1)}\ar[ddr] &\ar[u] W_{2,q^{-3}}^{(2)} \ar[rd]&
\\
&\ar[u] \ar[rd] W_{3,q^{-5}}^{(2)} & \fbox{$\scriptstyle W_{1,q^{-2}}^{(1)}$} 
\\
  &\fbox{$\scriptstyle W_{4,q^{-7}}^{(2)}$} \ar[u]\ar[dl]& \fbox{$\scriptstyle W_{1,q^{-4}}^{(1)}$} 
\\
\fbox{$\scriptstyle W_{2,q^{-6}}^{(1)}$} \ar[uuu]&&
}
\]
\caption{\label{Figg} {\it The isomorphism $\iota_\ell$ in type $G_2$ for $\ell = 3$.}}
\end{figure}

\begin{example}
{\rm
We take $\g$ of type $G_2$, and we choose $\ell = 3$.
The isomorphism $\iota_3$ is 
displayed in Figure~\ref{Figg}, which shows the image of the initial seed of $\AA_3$. 
There are $4$ frozen variables. $\AA_3$ has finite cluster type $A_4$. So $\CC_3$ has 18 prime objects, namely, the 8 Kirillov-Reshetikhin
modules of the initial seed, 
%
together with
\[
W_{1,q^{-6}}^{(1)},\ 
W_{1,q^{-3}}^{(2)},\
W_{1,q^{-5}}^{(2)},\ 
W_{2,q^{-5}}^{(2)},\ 
W_{2,q^{-7}}^{(2)},\ 
W_{3,q^{-7}}^{(2)},\ 
W_{1,q^{-7}}^{(2)},\
L((\varpi_1,q^0)+(\varpi_2,q^{-7})),\]
\[L((\varpi_1,q^0)+(\varpi_2,q^{-5})+(\varpi_2,q^{-7})),\
L((\varpi_1,q^0)+(\varpi_2,q^{-3})+(\varpi_2,q^{-5})+(\varpi_2,q^{-7})).
\]
}
\end{example}

\begin{remark} 
{\rm
(i)\ There is a similar Conjecture in which the finite quiver $\Gamma_\ell$ is replaced by the semi-infinite quiver $\Gamma^-$ 
with vertex set 
\[
V^- :=\{(i,r)\in V|r + d_i\leq 0\},
\]
and the category $\mathcal{C}_\ell$ is replaced by the category $\mathcal{C}^-$ of finite-dimensional $U_q(\hg)$-modules whose 
composition factors are of the form $L(\hla)$ with $\hla\in\bigoplus_{(i,r)\in V^-}\N (\varpi_i,q^{r+d_i})$, see \cite{HL2}.

\medskip\noindent
(ii)\ There is also a similar Conjecture in which the finite quiver $\Gamma_\ell$ is replaced by the doubly-infinite quiver $\Gamma$. 
In that case, the corresponding category is no longer a subcategory of the category of finite-dimensional $U_q(\hg)$-modules.
We have to consider a certain subcategory of the category $\mathcal{O}$ of (possibily infinite-dimensional) representations over a quantum Borel
subalgebra $U_q(\widehat{\mathfrak{b}})$ of $U_q(\hg)$, see \cite{HL3}. 
The category $\CC_\Z$ of Definition~\ref{defCZ} 
can be regarded as a subcategory of this category of $U_q(\widehat{\mathfrak{b}})$-modules.
The initial seed consists of the classes of prefundamental representations, which are simple infinite-dimensional modules 
of $U_q(\widehat{\mathfrak{b}})$ which cannot be extended to $U_q(\hg)$-modules.
}
\end{remark}

\subsection{What is known?}

\subsubsection{Part (ii)}
The difficult part of Conjecture~\ref{mainConj} is (i). 
If (i) is known, then (ii) follows from a result of \cite{GLS} which says that if a cluster algebra 
is a factorial ring, then every cluster variable is a prime element of this ring.

\subsubsection{First evidences}
As explained in \S\ref{SectSl2}, when $\g = \Sl_2$ Conjecture~\ref{mainConj} follows from \cite{CP1}.

In \cite{HL1}, (i) was proved for type $A_n$ and $D_4$ when $\ell = 1$. The proof was algebraic and combinatorial,
and certain parts of the proof were more general. 

In \cite{N1}, Nakajima proved (i) for types $A, D, E$ and $\ell = 1$, using the geometric approach to $U_q(\hg)$
via quiver varieties. (An introduction to this proof is presented in \cite{L}.)

A variant of Conjecture~\ref{mainConj} for type $A_n$ and $D_n$ when $\ell = 1$
was proved in \cite{HLJimboProc}. This involves finite subquivers of $\Gamma$ different from $\Gamma_\ell$.
Very recently Brito and Chari \cite{BC} generalized the results of \cite{HLJimboProc} in type $A$ using 
purely representation theoretic methods.

\subsubsection{Proof in simply-laced cases}
In \cite{Q}, Qin gave a proof of (i) for types $A, D, E$ and arbitrary $\ell$.
The proof also relies on the geometric approach, and uses the $t$-deformation of $K_0(\CC_\ell)$ introduced 
by Varagnolo-Vasserot and Nakajima in terms of quiver varieties.

\subsubsection{Connection with quiver Hecke algebras}\label{quiverH}
In type $A,D,E$, for $\ell = h/2 -1$, where $h$ is the Coxeter number supposed to be even, there is another
approach as follows. 
In \cite{HLCrelle}, we have shown that there is a ring isomorphism
\[
 i: \C\otimes_\Z K_0(\CC_{h/2-1}) \longrightarrow \C[N],
\]
where $N$ is a maximal unipotent subgroup of a simple Lie group $G$ with $\mathrm{Lie}(G) = \g$.
The ring of polynomial functions $\C[N]$ has a well-known cluster algebra structure, and the 
isomorphism $i$ transports the cluster structure of $\AA_{h/2-1} \cong \C\otimes_\Z K_0(\CC_{h/2-1})$
to the cluster structure of $\C[N]$. We have shown that the isomorphism $i$ maps the basis of 
$\C\otimes_\Z K_0(\CC_{h/2-1})$ consisting of classes of simple objects to the dual canonical basis
(or upper global basis) of $\C[N]$. Therefore to prove (i) in this case amounts to prove:
\begin{itemize}
 \item[(i')] The cluster monomials of $\C[N]$ form a subset of the dual canonical basis of $\C[N]$.
\end{itemize}
This was proved by Kang, Kashiwara, Kim and Oh \cite{KKKO}. They used the categorification of the dual
canonical basis of $\C[N]$ by simple objects of a category $\H$ of graded modules over quiver Hecke algebras.
This raises the question of a relation between the two categories $\H$ and $\CC_{h/2-1}$.
In \cite{KKK}, Kang, Kashiwara and Kim constructed a functor from $\H$ to $\CC_{h/2-1}$ inducing
the isomorphism $i^{-1}$ at the level of Grothendieck rings. In type $A$ this can be regarded as 
a variant of the quantum affine Schur-Weyl duality of Chari-Pressley and Ginzburg-Reshetikhin-Vasserot \cite{CPSchur,GRV}. 
Recently, Fujita \cite{F1,F2} proved that the KKK-functor is in fact an equivalence of categories.

\subsubsection{Non simply-laced cases}

\begin{figure}[t]
\[
\def\objectstyle{\scriptstyle}
\def\lablestyle{\scriptstyle}
\xymatrix@-1.0pc{
&&\ar[ld] (2,-1) \ar[rd]&
\\
&(1,-2) \ar[rd] &&\ar[ld] (3,-2)
\\
&&\ar[uu]\fbox{$\scriptstyle(2,-3)$} &
\\
&\ar[uu]\fbox{$\scriptstyle(1,-4)$} && \fbox{$\scriptstyle(3,-4)$}\ar[uu]
} 
\xymatrix@-1.0pc{
&&\ar[ld] (2,-1) &
\\
&(1,-2) \ar[rd] &&\ar[lu] (3,-2)\ar[dd]
\\
&& \fbox{$\scriptstyle(2,-3)$}\ar[ur] &
\\
&\ar[uu]\fbox{$\scriptstyle(1,-4)$} && \fbox{$\scriptstyle(3,-4)$}
}
\xymatrix@-1.0pc{
&&\ar[rd] (2,-1) &
\\
&(1,-2) \ar[ru]\ar[rd] &&\ar[ll] (3,-2)\ar[dd]
\\
&& \fbox{$\scriptstyle(2,-3)$}\ar[ur] &
\\
&\ar[uu]\fbox{$\scriptstyle(1,-4)$} && \fbox{$\scriptstyle(3,-4)$}
} 
\]
\[\def\objectstyle{\scriptstyle}
\def\lablestyle{\scriptstyle}
\xymatrix@-1.0pc{
&&\ar[ld] (2,-1) &&
\\
&(1,-2) \ar[rr]\ar[dd] && (3,-2)\ar[dd]&
\\
&& \fbox{$\scriptstyle(2,-3)$}\ar[ul] &&=
\\
&\ar[uuur]\fbox{$\scriptstyle(1,-4)$} && \fbox{$\scriptstyle(3,-4)$}&
}
\xymatrix@-1.0pc{
&&&\\
&&\ar[ld] (3,-2) &\ar[ld] (2,-1)  
\\
&\fbox{$\scriptstyle(3,-4)$} & \ar[u]\ar[rd](1,-2)&
\\
&&\ar[u] \fbox{$\scriptstyle(2,-3)$} & \fbox{$\scriptstyle(1,-4)$} \ar[uu]  
}
\]

\caption{\label{Fig7} {\it Mutations : from $\mathcal{A}_1$ in type $A_3$ to $\mathcal{A}_2$ in type $B_2$.}}
\end{figure}
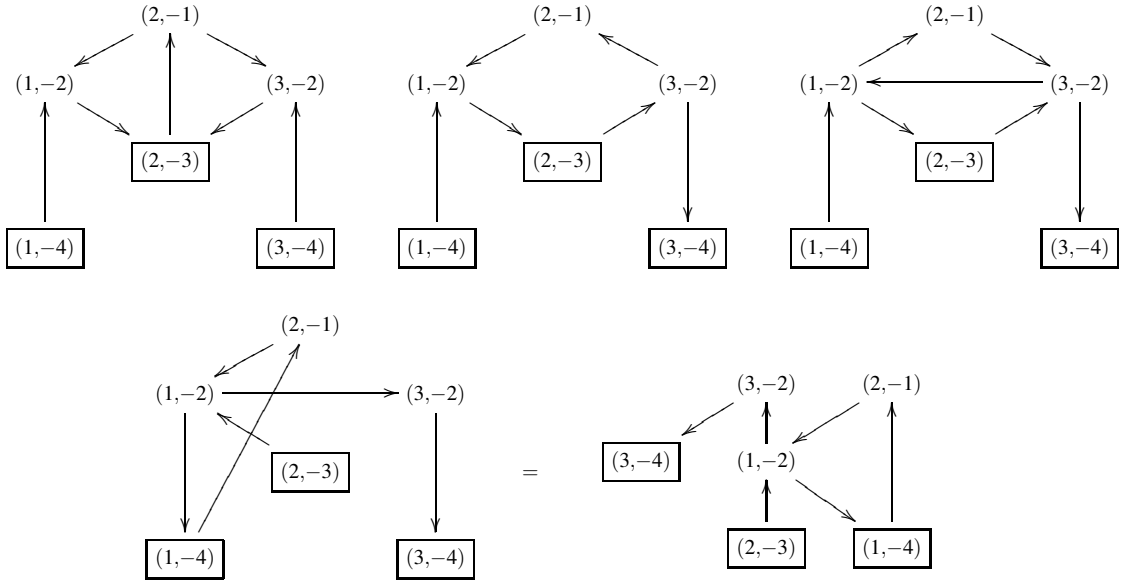

For non simply-laced types, let us start with the example of the category $\mathcal{C}_2$ in type $B_2$ 
discussed in Example \ref{example_conjecture_B2_l1}.
Comparing with the category $\mathcal{C}_1$ in type $A_3$ discussed in Example \ref{example_conjecture_A3_l1}, we observe that 
not only the cluster algebras $\mathcal{A}_1$ in type $A_3$ and $\mathcal{A}_2$ in type $B_2$ have the same cluster type $A_3$ and the
same numbers of frozen variables, but also the quivers in Figure~\ref{Fig2} and Figure~\ref{Fig3} are mutation-equivalent. 
This is illustrated in Figure~\ref{Fig7} with the mutation sequence at nodes $(3,-2)$,  $(2,-1)$, $(1,-2)$. 
Arrows between frozen vertices may be omitted 
as this does not change the cluster algebra structure. Hence we get a distinguished isomorphism between Grothendieck rings 
$$K_0(\mathcal{C}_{1,\, A_3}) \simeq K_0(\mathcal{C}_{2,\,B_2})$$ 
which is compatible with the cluster algebra structures. 
This last point is the most important since we know already that the two rings are isomorphic to the polynomial ring in $6$ variables.

This example is the first instance of a family of isomorphisms of cluster algebras 
\[
K_0(\mathcal{C}_{h/2 - 1,\, A_{2n-1}}) \simeq K_0(\mathcal{C}'_{B_n}) 
\]
obtained by the first author in a joint work with Hironori Oya \cite{HO}. 
These are distinguished isomorphisms between the Grothendieck rings of the type $A_{2n-1}$ categories $\mathcal{C}_{h/2-1,\, A_{2n-1}}$  
already mentioned in Section \ref{quiverH} above, 
and remarkable subcategories $\mathcal{C}'_{B_n}$ of finite-dimensional representations in
type $B_n$. The proof is established at the level of quantum cluster algebras, in order
to demonstrate a conjecture on quantum Grothendieck rings and related analogs of Kazhdan-Lusztig polynomials formulated in \cite{Hqgr}.

Using a completely different method based on functors from categories of representations of quiver Hecke algebras, 
Kashiwara-Kim-Oh and Kashiwara-Oh \cite{KKO,KO} constructed isomorphisms of Grothendieck rings in types $A_{2n-1}/B_n$ preserving the classes of simple modules. 
In fact, these match the distinguished isomorphisms obtained from cluster algebra structures in \cite{HO}. Hence the cluster algebra isomorphisms
also preserve classes of simple modules. For instance, in the example above, the bijection between prime simples can be directly written 
from the bijection with almost positive roots in Examples \ref{example_conjecture_A3_l1} and \ref{example_conjecture_B2_l1}. 
Consequently, combining these results, the analogue of Conjecture \ref{mainConj} for the subcategories $\mathcal{C}'_{B_n}$ holds.

It is expected that this approach will be extended to larger categories and to more general types.

\subsubsection{Real modules}
When the cluster algebra $\AA_\ell$ is not of finite cluster type, cluster monomials do not span
the vector space $\C \otimes \AA_{\ell}$. This raises the question of describing the simple objects
of $\CC_\ell$ whose class in the Grothendieck ring is a cluster monomial.

\begin{Conj}{\rm\cite{HL1, HL2}}\label{ConjReal}
The class of a simple object $S$ of $\CC_\ell$ in the Grothendieck
ring is a cluster monomial if and only if S is a \emph{real simple object},
that is, if and only if $S\otimes S$ is simple.
\end{Conj}

Let us assume that Conjecture~\ref{mainConj}~(i) holds. Then
one direction of Conjecture~\ref{ConjReal} is obvious: the square of a cluster monomial is clearly a cluster monomial.
The converse is wide open.

Real $U_q(\hg)$-modules have interesting properties. For instance Kang, Kashiwara, Kim and Oh
proved the following theorem, which was conjectured in type $A$ in \cite{LeImag}.

\begin{Thm}{\rm\cite{KKKOComp}}\label{ThmReal}
If $S_1$ and $S_2$ are two simple $U_q(\hg)$-modules, and one of them (at least) is real, then
$S_1\otimes S_2$ has a simple socle and a simple head. Moreover the socle and the head are isomorphic
if and only if $S_1\otimes S_2$ is simple.
\end{Thm}

Classification of real simple modules (in terms of their highest loop-weight) is a difficult open problem.
Recently Lapid and Minguez \cite{LM} classified in type $A$ all real simples satisfying a certain
regularity condition. Surprisingly, this classification is related to the classification of 
rationally smooth Schubert varieties in type $A$ flag varieties.


\section{Geometric character formulas}

An important obstacle for proving Conjecture~\ref{mainConj} in general is the absence of
Nakajima's geometric theory in the non-symmetric cases $B_n$, $C_n$, $F_4$, $G_2$.
It turns out that, applying the results of Derksen, Weyman, and Zelevinsky
\cite{DWZ1,DWZ2} to the cluster algebras $\AA_\ell$, 
one can define projective varieties whose Euler characteristics calculate the
$q$-characters of the \emph{standard} $U_q(\hg)$-modules in all types. These varieties can be seen
as generalizations of the Nakajima graded varieties $\LL^\bullet(V,W)$ for types $A, D, E$.
We shall now review this theorem of \cite{HL2}.

\subsection{Quiver Grassmannians and $F$-polynomials}

Let $Q$ be a quiver with vertex set $I$. Let $M$ be a representation of $Q$ over the field $\C$ of complex numbers.
Let $\be=(e_i)_{i\in I} \in \N^I$ be a dimension vector, and write $e := \sum_{i\in I} e_i$.
The variety $\Gr(\be,M)$ is the closed subvariety of the Grassmannian $\Gr(e,M)$ of $e$-dimensional
subspaces of $M$ whose points parametrize the sub-representations of $M$ with dimension vector $\be$. Thus, 
$\Gr(\be,M)$ is a projective complex variety, called a \emph{quiver Grassmannian}.

\begin{Def}
 The \emph{$F$-polynomial} of the representation $M$ of $Q$ is 
 \[
  F_M(\bv) := \sum_{\be \in \N^I} \chi(\Gr(\be,M)) \bv^\be,
 \]
where $\bv := (v_i)_{i\in I}$ is a sequence of commutative variables, $\bv^\be := \prod_i v_i^{e_i}$,
and $\chi(V)$ denotes the Euler characteristic of a complex projective variety $V$.
\end{Def}

\subsection{The algebra $A$}

Recall that in \S\ref{ssectQuivers} we have associated with $U_q(\hg)$ an infinite quiver $\G$.
Recall also the notation of \S\ref{ssecCartan} for the Cartan matrix. For every negative entry
$c_{ij}<0$ of the Cartan matrix and every $(i,r)\in V$, the graph $\G$ contains an oriented cycle
$\gamma_{i,j,r}$:
\begin{equation}\label{cycle}
\xymatrix@-1.0pc{
(i,r)\ar[rrdd] \\
(i,r-b_{ii})\ar[u]\\
{}\save[]{\vdots}\restore&&(j,r+b_{ij})\ar[lldd]\\
(i,r+2b_{ij}+b_{ii})\\
(i,r+2b_{ij})\ar[u]
}
\end{equation}

We define a \emph{potential} $S$ as the formal sum of all these
oriented cycles $\gamma_{i,j,r}$ up to cyclic permutations, see \cite[\S3]{DWZ1}. 
This is an infinite sum, but note that a given arrow
of $\G$ can only occur in a finite number of summands.
Hence all the cyclic derivatives of $S$, defined as in \cite[Definition 3.1]{DWZ1}, 
are finite sums of paths in~$\G$. 
Let $R$ be the list of all cyclic derivatives of $S$.
Let $J$ denote the two-sided
ideal of the path algebra $\C\G$ generated by $R$.
Following \cite{DWZ1}, we now introduce 
\begin{Def}
Let $A$ be the infinite-dimensional $\C$-algebra $\C\G\slash J$. 
\end{Def}

\begin{figure}
\[
\def\objectstyle{\scriptstyle}
\xymatrix@-1.0pc{
{}\save[]-<0cm,2ex>*{\vdots}\restore&\ar[ld]^{\a}(2,5)\\
{(1,4)}\ar[rd]^{\beta}& 
\\
&\ar[ld]^{\a} (2,3)\ar[uu]^{\gamma} 
\\
{(1,2)}\ar[rd]^{\beta}\ar[uu]^{\delta}& 
\\
&\ar[ld]^{\a} (2,1) \ar[uu]^{\gamma}
\\
\ar[uu]^{\delta}{(1,0)}\ar[rd]^{\beta}&
\\
&\ar[ld]^{\a}\ar[uu]^{\gamma} (2,-1) 
\\
\ar[uu]^{\delta}{(1,-2)}
&{}\save[]+<0cm,2ex>*{\vdots}\restore
}
\]
\caption{\label{Fig4} {\it The algebra $A$ in type $A_2$.}}
\end{figure}
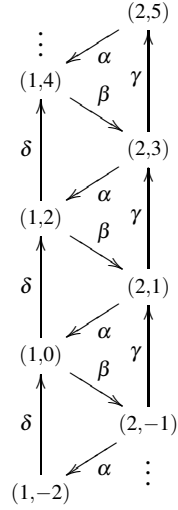

\begin{example}\label{example_algebra_A}
{\rm
Let $\g=\Sl_3$, of type $A_2$. 
The quiver $\G$ is displayed in Figure~\ref{Fig4}. We slightly abuse notation using the same letter $\a$ for every arrow of the form
$(2,r) \to (1,r-1)$, and similarly for $\b$, $\ga$, $\de$. All cycles $\ga_{i,j,r}$ are of length 3, and are either of the form 
$(\ga, \b, \a)$ or of the form $(\b, \de, \a)$. Therefore the potential is
\[
 S = \sum (\ga, \b, \a) +  \sum (\b, \de, \a),
\]
and the relations obtained by taking its cyclic derivatives are of the form:
\[
\begin{array}{lcc}
\mbox{derivative with respect to $\a$}& \leadsto & \ga\b + \b\de = 0, \\[3mm]
\mbox{derivative with respect to $\b$}& \leadsto & \a\ga + \de\a = 0, \\[3mm]
\mbox{derivative with respect to $\ga$}& \leadsto & \b\a = 0, \\[3mm]
\mbox{derivative with respect to $\de$}& \leadsto & \a\b  = 0.
\end{array} 
\]
So $A$ is the algebra defined by the quiver $\G$ of Figure~\ref{Fig4}, subject to the above 4 families of relations.
}
\end{example}

\subsection{Some $A$-modules}

The algebra $A$ is infinite-dimensional. For every $(i,r)\in V$ there is a one-dimensional $A$-module $S_{(i,r)}$ supported
on vertex $(i,r)$. Let $I_{(i,r)}$ denote the injective hull of $S_{(i,r)}$, an infinite-dimensional indecomposable $A$-module.

\begin{figure}
\[
\def\objectstyle{\scriptstyle}
\xymatrix@-1.0pc{
\\
&\ar[ld]^{\a} (2,3) 
\\
{(1,2)}& 
\\
&\ar[ld]^{\a} (2,1) \ar[uu]^{\gamma}
\\
\ar[uu]^{\delta}{(1,0)}&
\\
&\ar[ld]^{\a}\ar[uu]^{\gamma} (2,-1) 
\\
\ar[uu]^{\delta}{(1,-2)}
&{}\save[]+<0cm,2ex>*{\vdots}\restore
\\
{}\save[]+<0cm,2ex>*{\vdots}\restore&
}
\]
\caption{\label{Fig5} {\it The injective $A$-module $I_{(1,2)}$ in type $A_2$.}}
\end{figure}
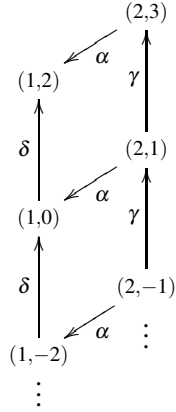

\begin{example}\label{example_algebra_A_cont}
{\rm
We continue Example~\ref{example_algebra_A}. The injective module $I_{(1,2)}$ is represented in 
Figure~\ref{Fig5}. Each vertex occuring in the picture carries a one-dimensional vector space, 
and all occuring arrows are nonzero. The $\beta$ arrows are all zero and therefore not represented.
The infinite socle series of this module is:
\[
\begin{array}{c}
 \vdots \\[3mm]
 S_{(1,-2)}\oplus S_{(2,1)} \\[3mm]
 S_{(1,0)}\oplus S_{(2,3)} \\[3mm]
 S_{(1,2)}
\end{array}
\]
Note that we have a short exact sequence of $A$-modules
\[
0 \to K_{(1,2)} \to I_{(1,2)} \to I_{(1,0)} \to 0, 
\]
where $K_{(1,2)}$ is the two-dimensional module:
\[
\def\objectstyle{\scriptstyle}
\xymatrix@-1.0pc{
&\ar[ld]^{\a} (2,3) 
\\
{(1,2)}& 
}
\]
}
\end{example}

\begin{Prop}
For every $(i,r)\in V$ there is a unique submodule $K_{(i,r)}$ of $I_{(i,r) }$ such that $I_{(i,r)}/K_{(i,r)}$ is
isomorphic to $I_{(i,r-d_i)}$. The module $K_{(i,r)}$ is finite-dimensional. 
\end{Prop}

\subsection{Geometric $q$-character formula}
The next theorem gives a geometric description of the $q$-character of the fundamental $U_q(\hg)$-module $L((\varpi_i,q^{r-d_i}))$
in terms of the $F$-polynomial of the $A$-module $K_{(i,r)}$. To state it we need some more notation. 
To every $(i,r)\in V$ we attach a commutative variable $z_{(i,r)}$, and we set
\[
 \hy_{(i,r)} := \prod_{(i,r)\to (j,s)} z_{(j,s)} \prod_{(k,l)\to(i,r)}z_{(k,l)}^{-1},
\]
where the first product runs over all arrows of $\G$ starting at vertex $(i,r)$, and the second product over
all arrows of $\G$ ending in vertex $(i,r)$.

\begin{Thm}{\rm\cite{HL2}}\label{ThmGeo}
The $F$-polynomial $F_{K_{(i,r)}}(\hby)$ of the $A$-module $K_{(i,r)}$ evaluated in the variables $\hy_{(i,r)}$ can be expressed 
as a Laurent polynomial in the new variables
\[
 Y_{j,q^{s-d_j}} := \frac{z_{(j,s-2d_j)}}{z_{(j,s)}},\qquad ((j,s)\in V).
\]
Then we have
\[
\chi_q(L((\varpi_i,q^{r-d_i}))) = Y_{i,q^{r-d_i}} F_{K_{(i,r)}}(\hby) 
\]
where in the right-hand side $F_{K_{(i,r)}}(\hby)$ is expressed in terms of the variables $Y_{j,q^{s-d_j}}$.
\end{Thm}

\begin{example}
{\rm
We continue Example~\ref{example_algebra_A_cont}.
The $A$-module $K_{(1,2)}$ has exactly three submodules: $\{0\}$, $S_{(1,2)}$, and $K_{(1,2)}$. 
So there are three nonempty quiver Grassmannians, each reduced to a single point, hence each having Euler characteristic 1.
It follows that
\[
 F_{K_{(1,2)}}(\bv) = 1 + v_{(1,2)} + v_{(1,2)}v_{(2,3)}.
\]
On the other hand
\[
 \hy_{(1,2)} = \frac{z_{(2,1)}z_{(1,4)}}{z_{(1,0)}z_{(2,3)}},\qquad
 \hy_{(2,3)} = \frac{z_{(1,2)}z_{(2,5)}}{z_{(1,4)}z_{(2,1)}}.
\]
Hence
\[
 F_{K_{(1,2)}}(\by) = 1 + \frac{z_{(2,1)}z_{(1,4)}}{z_{(1,0)}z_{(2,3)}} + \frac{z_{(1,2)}z_{(2,5)}}{z_{(1,0)}z_{(2,3)}}
 = 1 + Y_{1,q}^{-1}Y_{1,q^3}^{-1}Y_{2,q^2} +  Y_{1,q}^{-1} Y_{2,q^4}^{-1},
\]
and
\[
 Y_{1,q}F_{K_{(1,2)}}(\by) = Y_{1,q} + Y_{1,q^3}^{-1}Y_{2,q^2} + Y_{2,q^4}^{-1} = \chi_q(L((\varpi_1,q))).
\]
}
\end{example}

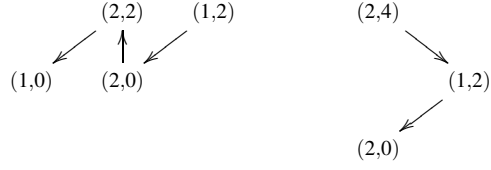
\begin{figure}
\[
\def\objectstyle{\scriptstyle}
\xymatrix@-1.0pc{
&&\ar[ld] (\bdeux,2)& \ar[ld](\bun,2)&
\\
&{(\bun,0)} &\ar[u] (\bdeux,0)&&
}
\qquad
\xymatrix@-1.0pc{
(\bdeux,4) \ar[rd]&&
\\
 &\ar[ld] (\bun,2)   &
\\
(\bdeux,0) &&
}
\]
\caption{\label{Fig6} {\it The $A$-modules $K_{(1,0)}$ and $K_{(2,0)}$ in type $B_2$.}}
\end{figure}

\begin{example}
{\rm
We take $\g = \mathfrak{so}_5$, of type $B_2$.
The $A$-modules $K_{(1,0)}$ and $K_{(2,0)}$ are displayed in Figure~\ref{Fig6}.
In this case too, the nonempty quiver Grassmannians are reduced to points, and we can calculate
easily:
\[
\begin{array}{lcl}
\chi_q(L((\varpi_1,q^{-2}))& = & Y_{1,q^{-2}}\left(1 + \hy_{(1,0)} + \hy_{(1,0)}\hy_{(2,2)} + \hy_{(1,0)}\hy_{(2,2)}\hy_{(2,0)}
+ \hy_{(1,0)}\hy_{(2,2)}\hy_{(2,0)}\hy_{(1,2)} \right),\\[2mm]
\chi_q(L((\varpi_2,q^{-1}))& = & Y_{2,q^{-1}}\left(1 + \hy_{(2,0)} + \hy_{(2,0)}\hy_{(1,2)} + \hy_{(2,0)}\hy_{(1,2)}\hy_{(2,4)}
\right).
\end{array}
\]
To express these $q$-characters in terms of variables $Y_{j,q^s}$ one uses the formulas:
\[
\hy_{(1,r)} = \frac{z_{(1,r+4)}z_{(2,r-2)}}{z_{(1,r-4)}z_{(2,r+2)}} = \frac{Y_{2,q^{r-1}}Y_{2,q^{r+1}}}{Y_{1,q^{r+2}}Y_{1,q^{r-2}}},
\qquad
\hy_{(2,r)} = \frac{z_{(1,r-2)}z_{(2,r+2)}}{z_{(1,r+2)}z_{(2,r-2)}} = \frac{Y_{1,q^{r}}}{Y_{2,q^{r-1}}Y_{2,q^{r+1}}}.
\]
}
\end{example}

\subsection{Comments on Theorem~\ref{ThmGeo}}

\subsubsection{Kirillov-Reshetikhin modules}
In \cite{HL2} we give a similar $q$-character formula for every Kirillov-Reshetikin module.
One only needs to replace the $A$-modules $K_{(1,r)}$ by some larger finite-dimensional submodules of the injective
modules $I_{(i,r)}$.

\subsubsection{Standard modules}
Classical properties of Euler characteristics imply that, given two finite-dimensional $A$-modules $M$ and $N$, we have
\[
 F_{M\oplus N}(\bv) = F_M(\bv) F_N(\bv).
\]
On the other hand, $q$-characters are multiplicative on tensor products. So Theorem~\ref{ThmGeo} readily extends to
tensor products of fundamental $U_q(\hg)$-modules, that is, we have a similar geometric $q$-character formula for
\emph{standard} $U_q(\hg)$-modules (or local Weyl modules), in which one uses quiver Grassmannians of direct sums of $A$-modules $K_{(i,r)}$.

\subsubsection{Relation to Nakajima's theory}
If $\g$ is of type $A,D,E$ it follows from results of Lusztig \cite{Lu} and Savage-Tingley \cite{ST} that
the quiver Grassmannians
\[
 G(\be,\ba) := \Gr\left(\be,\bigoplus_{(i,r)\in V}K_{(i,r)}^{\oplus a_{(i,r)}}\right)
\]
are homeomorphic to Nakajima varieties $ \LL^\bullet(V,W)$, 
where the graded dimension of $V$ is encoded in the dimension vector $\be$, and
the graded dimension of $W$ is given by the multiplicity vector $\ba = (a_{(i,r)})$.
One can therefore regard the varieties $G(\be,\ba)$ in non simply-laced types $B, C, F, G$ as 
natural candidates for replacing the graded Nakajima varieties $ \LL^\bullet(V,W)$.

\subsubsection{Beyond KR-modules and standard modules}
By the Derksen-Weyman-Zelevinsky theory, every cluster monomial of $\AA_\ell$ has an expression
of the form
\[
 m = \bz^\bg \,F_K(\hby)
\]
for an appropriate $A$-module $K$.
So, if Conjecture~\ref{mainConj} is true, all the simple $U_q(\hg)$-modules corresponding to
cluster monomials (all the real modules, if Conjecture~\ref{ConjReal} is true) have a similar geometric $q$-character formula in terms
of quiver Grassmannians.




\bigskip
\small
\noindent
\begin{tabular}{ll}
David {\sc Hernandez}  
& Sorbonne Paris Cit\'e, Univ Paris Diderot,\\
&CNRS Institut de
  Math\'ematiques de Jussieu-Paris Rive Gauche UMR 7586,\\
& B\^atiment Sophie Germain, Case 7012,
75205 Paris Cedex 13, France\\
&Institut Universitaire de France,\\
& email : {\tt david.hernandez@imj-prg.fr}\\
[5mm]
Bernard {\sc Leclerc}  & Universit\'e de Caen Normandie,\\
&CNRS UMR 6139 LMNO, 14032 Caen, France\\
&email : {\tt bernard.leclerc@unicaen.fr}
\end{tabular}

\end{document}